# Towards Homological Methods in Graphic Statics

Zoe COOPERBAND[*,a] and Robert GHRIST[a,b]

[*]Department of Electrical and Systems Engineering
University of Pennsylvania, 200 S. 33rd St., Philadelphia, PA, USA
{zcoop,ghrist}@seas.upenn.edu

[a] Department of Electrical and Systems Engineering, University of Pennsylvania
[b] Department of Mathematics, University of Pennsylvania

**Abstract**
Recent developments in applied algebraic topology can simplify and extend results in graphic statics – the analysis of equilibrium forces, dual diagrams, and more. The techniques introduced here are inspired by recent developments in cellular cosheaves and their homology. While the general theory has a few technical prerequisites (including homology and exact sequences), an elementary introduction based on little more than linear algebra is possible. A few classical results, such as Maxwell's Rule and 2D graphic statics duality, are quickly derived from core ideas in algebraic topology. Contributions include: (1) a reformulation of statics and planar graphic statics in terms of cosheaves and their homology; (2) a new proof of Maxwell's Rule in arbitrary dimensions using Euler characteristic; and (3) derivation of a novel relationship between mechanisms of the form diagram and obstructions to the generation of force diagrams. This last contribution presages deeper results beyond planar graphic statics.

## 1. Introduction

*Applied topology* is a branch of contemporary applied mathematics that uses algebraic topology as a toolset in applications [1]. Some prominent examples from the past two decades include the use of topological tools in data analysis [2], neuroscience [3], sensor networks [4], and signal processing [5].

One of the more recent and exciting developments in applied topology is the establishment and implementation of *cellular sheaf theory*. Building on the classical and powerful sheaf theory of the 1950s-60s [6], cellular sheaves (and dual *cosheaves*) provide a locally finite and inherently computation-friendly toolset which matured in the thesis of Curry [7]. These tools have been used in network coding [8], pursuit and evasion games [9], distributed optimization [10], opinion dynamics [11], scene analysis [12], and extending spectral graph theory [13].

In this introduction we use cellular cosheaves to characterize the static behavior of trusses. Following a classical history [14], [15], advances in the late 20th century [16]–[18] laid foundations for homological methods in static analysis. Algebraic foundations for the mechanical properties of two dimensional graphic statics has been further developed since [19]–[21]. Additional great strides have been made in algebraic descriptions of 3D polyhedral graphic statics [22]–[24]. Our contributions via cellular cosheaves extend these threads.

## 2. Cellular cosheaves

In the simplest possible setting, cosheaves are algebraic structures attached to a *network* (or graph of vertices and edges) and comprised of *vector spaces* and *linear transformations*.

## 2.1. A general construction

Let $X = (V, E)$ denote a graph with vertex set $V$ and (undirected) edge set $E$, with each edge connecting distinct vertices. A *cosheaf* $\mathcal{F}$ on $X$ associates to each vertex $v \in V$ and each edge $e \in E$ a finite dimensional vector space $\mathcal{F}_v$ ($\mathcal{F}_e$ respectively). These vector spaces – called *stalks* – can be thought of as *local data*, and they can vary from vertex to vertex and edge to edge. The cosheaf consists not only of the underlying graph along with the stalks, but also linear transformations connecting the stalks where edges and vertices meet. For each incident edge-vertex pair $e \triangleright v$, the cosheaf contains a linear transformation $\mathcal{F}_{e \triangleright v}: \mathcal{F}_e \to \mathcal{F}_v$ from data over $e$ to data over $v$. In all examples in this introduction, the vector spaces have explicit bases, and the maps $\mathcal{F}_{e \triangleright v}$ are matrices. One thinks of these *cosheaf matrices* as being programmable: they define the constraints of the cosheaf stalk data.

***Example 1:*** For any graph $X$, the unit *constant cosheaf* $\mathcal{I}_X$ has a 1-dimensional vector space as every stalk with all cosheaf matrices the identity. For a fixed (not necessarily 1-dimensional) vector space $V$, the $V$-constant cosheaf $V_X$ has a copy of $V$ on each vertex and edge, connected by identity matrices.

***Example 2:*** Consider single-variable *polynomial splines* over a graph. These splines over edges are degree $\leq m$ polynomials in a single variable, meeting at vertices with global smoothness class $C^r$ ($r$ times continuously differentiable). Consider the cosheaf $\mathcal{K}$ with edge stalks consisting of $P_m$, the vector space of degree $\leq m$ polynomials. The vertex stalks are quotient vector spaces $P_m/I^{r+1}$, where $I^{r+1}$ is the subspace of polynomials which vanish at the vertex with complementary degree at least $r + 1$. The matrices $\mathcal{K}_{e \triangleright v}$ are quotient projections $P_m \to P_m/I^{r+1}$. Further exposition may be found in [25].

Our definitions have focused on the case where the base space is an undirected network. Extension to (regular) higher dimensional cell complexes is straightforward [7], with stalks attached to higher dimensional cells and corresponding cosheaf matrices.

## 2.2 The force cosheaf

Our central object of study is the *force cosheaf* – a means of tracking force vectors across a truss. Fix a graph $X = (V, E)$ modelling a mechanical truss built from axial geometric members $E$ pinned at their ends with free-rotational joints $V$. Assume that $X$ is realized in Euclidean space $\mathbb{R}^n$ (with $n = 2,3$ being typical) by an assignment $p: V \to \mathbb{R}^n$ (that embeds edges linearly). The force cosheaf $\mathcal{F} = \mathcal{F}(X, p)$ for this realization of $X$ has vertex stalks $\mathcal{F}_v$ all equal to $\mathbb{R}^n$, encoding the static loads at the joints. The edge stalks $\mathcal{F}_e \cong \mathbb{R}$ encode a stress along the member: positive or negative for tension or compression respectively. For each edge $e$ connecting joints $u$ and $v$, the matrices $\mathcal{F}_{e \triangleright u}: \mathcal{F}_e \to \mathcal{F}_u$ and $\mathcal{F}_{e \triangleright v}: \mathcal{F}_e \to \mathcal{F}_v$ send the basis vector of $\mathcal{F}_e$ to the vectors $p(u) - p(v) \in \mathcal{F}_u$ and $p(v) - p(u) \in \mathcal{F}_v$ respectively. A sketch of the force cosheaf is pictured on the right in Figure 1.

The force cosheaf contains within it all the data from which statics and graphic statics unfolds. What reveals degrees of freedom and equilibrium is the classical algebraic-topological tool called *homology*.

## 3. Homology

Given a graph, one often considers its global topological properties such as the number of its connected components or independent cycles. These intuitive properties generalize greatly to homology – the algebraic study of *"holes."* Cosheaves, as introduced above, appear as mere containers of data; deeper meaning follows from considering their homology.

## 3.1 Classical cellular homology

We begin by describing a simple, classical instance of homology. In the case of $X = (V, E)$ a graph, one defines the vector spaces $C_0 = C_0 X$ and $C_1 = C_1 X$ with basis the formal set of vertices and edges respectively. Counting the number of edges and vertices we have $\dim C_0 = |V|$ and $\dim C_1 = |E|$. One may then define the *oriented incidence matrix* $\partial: C_1 \to C_0$ with respect to an arbitrary orientation given to edges. This matrix sends a basis edge $e$ to the linear combination of vertices $\partial e = u - v$, respecting this orientation. The null space, or *kernel* of this matrix $H_1 X = \ker \partial$ has basis of oriented cycles in $X$, where multiplication by $-1$ changes the orientation of the cycle. The quotient vector space $H_0 X = C_0 / \operatorname{im} \partial$ of $C_0$ by the image of the oriented incidence matrix is known as the *cokernel* (isomorphic to the left null space of $\partial$); it has basis of connected components of $X$.

This graph construction generalizes to $X$ a higher dimensional cell complex with more specialized language. The vector space $C_k = C_k X$ is known as the space of *k-chains*, linear combinations of formal $k$-cells (with arbitrary orientation). These vector spaces assemble into a sequence called a *chain complex* by means of *boundary matrices* $\partial: C_k \to C_{k-1}$ which send a basis $k$-cell $\sigma$ to the linear combination of boundary $(k-1)$-cells $\partial \sigma = \sum \pm \tau$, where the sign comes from compatibility of orientation [26]. Composed boundary matrices are the zero matrix, with $\partial_k \circ \partial_{k+1} = 0$.

The crucial definition is this. The (cellular) *homology* of $X$, $HX$, is a collection of vector spaces $H_k = H_k X$. This $k^{\text{th}}$ homology of $X$ is defined to be the quotient vector space $H_k X = \ker \partial_k / \operatorname{im} \partial_{k+1}$: kernels mod images, or equivalence classes of collections of $k$-cells with vanishing boundary, modulo distributions that are the boundary of a higher dimensional chain. The dimensions of the homologies, $\beta_k = \dim H_k X$, are known as the *Betti numbers* of $X$. In general, $\beta_0$ gives the number of connected components of $X$, and higher Betti numbers count the number of different dimensional holes within $X$ (in a graph, cycles are 1-dimensional holes). Homology is invariant with respect to the arbitrary orientation of cells and any homotopy (think: controlled deformation) of the space $X$.

### 3.2 Cosheaf homology

Homology extends to cosheaves. Let $\mathcal{F}$ be a cosheaf over a graph $X = (V, E)$. The space of cosheaf *0-chains* $C_0 = C_0 \mathcal{F} = \bigoplus_v \mathcal{F}_v$ is the product vector space bundling data over all vertices, and the space of *1-chains* $C_1 = C_1 \mathcal{F} = \bigoplus_e \mathcal{F}_e$ is the corresponding amalgamation of vector data over all edges. A vector in $C_0$ (or $C_1$) is interpreted as a distribution of data over vertices (or edges). The *cosheaf boundary matrix* $\partial: C_1 \to C_0$ acts on a 1-chain $\xi = (\xi_e)$ yielding a 0-chain $\partial \xi$ whose value on a vertex $v$ is given by

$$(\partial \xi)_v = \sum_{e_i \triangleright v} \pm \mathcal{F}_{e_i \triangleright v}(\xi_{e_i}) \tag{1}$$

where the sign is informed by an arbitrary choice of orientation and the sum is over all edges incident to $v$. Summations of the form (1) occur over all vertices, utilizing cosheaf matrices $\mathcal{F}_{e_i \triangleright v_j}$ to transform a distribution of data over edges to data over vertices. One imagines the vertices communicating and pulling data from their neighbor edges. For a cosheaf over a $n$-dimensional cell complex, the boundary matrix extends to a sequence of linear transformations $\partial_k: C_k \to C_{k-1}$, the entire sequence being collected into a chain complex $C\mathcal{F} = (C_k, \partial_k)$, with $\partial_k = 0$ when $k \leq 0$ or $k > n$.

The *homology* of the cosheaf, $H\mathcal{F}$, is the collection of vector spaces $H_k = H_k \mathcal{F}$, where $H_k \mathcal{F} = \ker \partial_k / \operatorname{im} \partial_{k+1}$. This is kernels mod images, just as in (classical) cellular homology. The difference between cellular homology and cosheaf homology is in the difference of data of the chain complex $CX$ or $C\mathcal{F}$ - the process of taking homology is afterwards identical.

*Examples:* Recall the examples in Section 2.1.

1. The cellular homology of a network $X$ is the same as the homology of the unit constant cosheaf $\mathcal{I}_X$ over it. Here, $H\mathcal{I}_X \cong HX$ and so the homology is the same as the cellular case. The Betti numbers of $\mathcal{I}_X$ count connected components and independent cycles of $X$. For a general constant cosheaf $V_X$, this has homology isomorphic to $m$ copies of $HX$, where $m = \dim V$.
2. The space of splines over a polygonal cycle is isomorphic to $H_1 \mathcal{K}$, the 1st homology of the spline cosheaf of Example 2. This is a simple example of more general results on splines [25].

These examples hint at the role that homology plays in classifying solutions of constrained systems.

### 3.3 Force cosheaf homology

What does the homology of the force cosheaf entail? At the physical level, each truss joint combines the axial forces along incident members. In the mathematical model, this is precisely what the force cosheaf boundary matrix describes. Given a 1-chain – a choice of axial force along all members – the boundary matrix $\partial: C_1\mathcal{F} \to C_0\mathcal{F}$, otherwise known as the *equilibrium matrix* [27], evaluates a distribution of net forces at all joints. This is significant: *the homology of the force cosheaf, $H\mathcal{F}$, characterizes equilibrium-stresses and configurations of the mechanism.* The details are worth emphasizing:

1. The 1-dimensional homology $H_1\mathcal{F} = \ker \partial$ consists of all stress distributions that keep the structure in static equilibrium (zero net force at each joint). The dimension of this null space (the 1st Betti number $\beta_1 = \dim H_1\mathcal{F}$) records the number of dimensions of *self-stress*.
2. The 0-dimensional homology, $H_0\mathcal{F} = C_0\mathcal{F}/\mathrm{im}\,\partial$, records distributions of forces at joints which cannot impart forces to edges. This space of un-resistible forces is interpretable as the space of constrained *degrees of freedom* of vertices of $X$, including global translations, rotations, and mechanisms. Under these infinitesimal motions, edge lengths would remain constant.

A sample graph $X$ with non-zero self-stress is provided in Figure 1. As is common in a well-constructed cosheaf, homology characterizes/classifies properties of the global system. With these tools, nontrivial inferences can be made in a unified fashion.

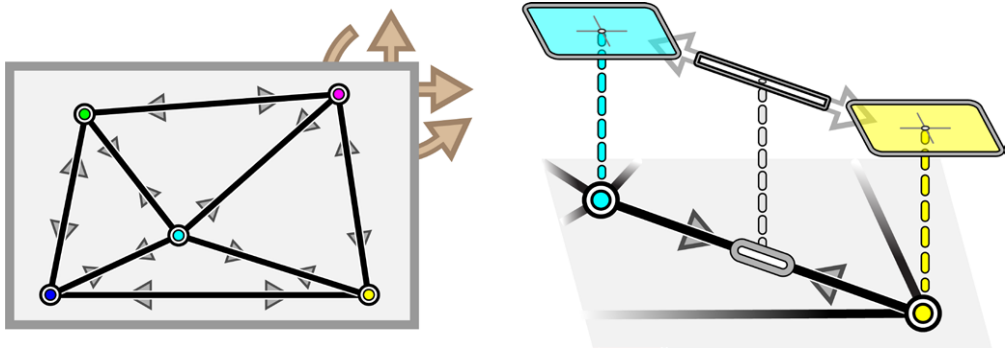

Figure 1. The force cosheaf $\mathcal{F}$ over the graph (left) has cosheaf boundary matrix of size 10 by 8 with nullity 1. The homology of $\mathcal{F}$ characterizes the self-stress over edges and three degrees of freedom (two translational and one rotational of the entire structure). A local view of $\mathcal{F}$ is highlighted (right), with the cosheaf matrices mapping edge stalks to vertex stalks "above" the graph.

## 4. Euler, Maxwell, & Betti

Maxwell's Rule [14] is the classic condition for a frame structure to be isostatic. The modern form for 3D trusses is given in [27]:

$$3|V| - |E| = 6 + |M| - |S|, \tag{2}$$

where $|M|$ is the number of linkage mechanisms, $|S|$ is the number of self-stresses, $|V|$ is the number of vertices of the truss, and $|E|$ is the number of its edges. The prevalence of differences (especially in the number of vertices and edges) is a clue that this is related to *Euler characteristic*.

Recall that elementary topological fact that for a finite polyhedral *spherical* surface $X$ with $|V|$ vertices, $|E|$ edges, and $|F|$ faces, the Euler characteristic $\chi(X) = |V| - |E| + |F|$ equals $+2$, independent of how the spherical surface is discretized. The reason for this is an elementary but fundamental result. For a finite collection of finite-dimensional vector spaces $C = \{C_k\}$, its Euler characteristic

$$\chi(C) = \sum_k (-1)^k \dim C_k \tag{3}$$

is likewise finite and agrees with the classical formula that counts vertices, edges, and faces. This is formulated as an alternating sum of Betti numbers by the following classical result:

**Lemma 1:** [26] *The Euler characteristic of a chain complex $C = (C_k)$ and its homology $H$ agree:*

$$\chi(C) = \chi(H) = \sum_{k}(-1)^k \beta_k, \tag{4}$$

This is why for a triangulated spherical surface $\chi = 2$, since its homology satisfies $\beta_0 = \beta_2 = 1$ and $\beta_1 = 0$. The proof of Lemma 1 is neither more nor less than the rank-nullity theorem from linear algebra, suitably generalized. Indeed, a single linear transformation $A: U \to V$ *is* an elementary chain complex where $C_1 = U$ and $C_0 = V$. Taking homology, one recovers the rank-nullity theorem as a special case.

We reinterpret Maxwell's Rule as coming from a chain complex and its homology: in this case, the force cosheaf $\mathcal{F}$. For a structure in $n$-D, its Euclidean group of translations and rotations in $n$-D has dimension $\binom{n+1}{2} = \frac{1}{2}n(n+1)$. With this, Lemma 1 applied to the force cosheaf $\mathcal{F}$ in $n$-D yields

$$n|V| - |E| = \dim C_0 - \dim C_1 = \dim H_0 - \dim H_1 = \binom{n+1}{2} + |M| - |S|. \tag{5}$$

*This is the Maxwell Rule in dimension $n$*, giving equation (2) when $n = 3$.

## 5. Boundary conditions

Trusses are outfitted with external loadings and reaction forces in practice. The statics of these trusses are characterized by *relative homology*, meaning homology relative to an external structure. This requires some additional theoretical tools.

### 5.1 Cosheaf maps and quotients

For two cosheaves $\mathcal{F}$ and $\mathcal{G}$ over a graph $X = (V, E)$, a *cosheaf map* $\phi: \mathcal{F} \to \mathcal{G}$ is comprised of stalk-wise linear maps $\phi_x: \mathcal{F}_x \to \mathcal{G}_x$ for $x$ ranging over every cell of $X$. In addition, the cosheaf map must respect the cosheaf matrices of $\mathcal{F}$ and $\mathcal{G}$: the equations $\phi_v \circ \mathcal{F}_{e \triangleright v} = \mathcal{G}_{e \triangleright v} \circ \phi_e$ must hold.

It may be the case that the data of cosheaf $\mathcal{F}$ can be considered as a subset of the data of $\mathcal{G}$. In this special case, $\phi: \mathcal{F} \hookrightarrow \mathcal{G}$ is an *inclusion* of the stalks of $\mathcal{F}$ into the stalks of $\mathcal{G}$, and one can think of $\mathcal{F}$ as lying "inside" $\mathcal{G}$. To model the data in $\mathcal{G}$ not captured by $\mathcal{F}$, one builds the *quotient cosheaf* $\mathcal{G}/\mathcal{F}$ with stalks of quotient vector spaces $(\mathcal{G}/\mathcal{F})_x = \mathcal{G}_x/\mathcal{F}_x$, with $x$ ranging over all cells. One can think of the cosheaf $\mathcal{G}/\mathcal{F}$ as encoding the *orthogonal compliment* to $\mathcal{F}$ inside of $\mathcal{G}$, satisfying the isomorphism of cosheaves $\mathcal{F} \oplus \mathcal{G}/\mathcal{F} \cong \mathcal{G}$ (comprised of vector space isomorphisms between stalks). Here, the data of $\mathcal{G}$ is subdivided into its constituent components parallel and orthogonal to $\mathcal{F}$. This construction will help with modelling external forces on a truss.

### 5.2 External loadings

We model external forces on a truss somewhat differently from other authors [19], [21], [27], considering boundary conditions as residing on an exterior truss abstracting the origin of the boundary loads. This exterior truss loosely set to be a loop around the structure: the intention of this is to model the loaded truss of interest as a sub-truss of a larger structure. The benefits of this approach are a clear differentiation between loaded vertices and unloaded vertices, as well as a homological formulation of equilibrium stresses. This is specifically important in simple 2D graphic statics, where loads are only applied on the exterior of a loaded truss.

Suppose that $Y = (V_Y, E_Y)$ (the exterior truss loop) is a closed subgraph of $X = (V, E)$ (the extended original truss with *boundary conditions*) as in Figure 2. For the force cosheaf $\mathcal{F} = \mathcal{F}_X$ with domain over $X$, the cosheaf $\mathcal{F}_Y$ has stalks identical to those of $\mathcal{F}_X$ for vertices and edges in $Y$ and zero stalks elsewhere.

It is straightforward to define the inclusion cosheaf map $\iota: \mathcal{F}_Y \hookrightarrow \mathcal{F}_X$, with stalk-wise linear maps the identity matrix over cells of $Y$ and the zero matrix over $X - Y$. The ensuing quotient cosheaf $\mathcal{F}_{X-Y} \coloneqq \mathcal{F}_X/\mathcal{F}_Y$ has support strictly on the set $X - Y$ (which is topologically an open disk). Some edges $\{e_l, e_r\}$ of $X - Y$ have an open end in $Y$ - these are *loaded edges* and *reaction edges* of $X - Y$. These edges represent lines of action of forces transferred between $Y$ and $X - Y$.

The quotient cosheaf $\mathcal{F}_{X-Y}$ models the static behavior of the original truss $X - Y$ with respect to the space of peripheral boundary conditions $Y$. The *relative homology* $H_1 \mathcal{F}_{X-Y}$ of this cosheaf characterizes the full space of *equilibrium stresses*, including internal self-stresses as well as external equilibrium loadings from edges $\{e_l, e_r\}$ (ignorant of whether $Y$ could deliver these forces or not). Here, the truss $X - Y$ is blind to the statics of $Y$ and liberally permits any possible loading through edges $\{e_l, e_r\}$. Our next goal is to relate the homologies of the force cosheaves $\mathcal{F}_X, \mathcal{F}_Y,$ and $\mathcal{F}_{X-Y}$ with varying domains.

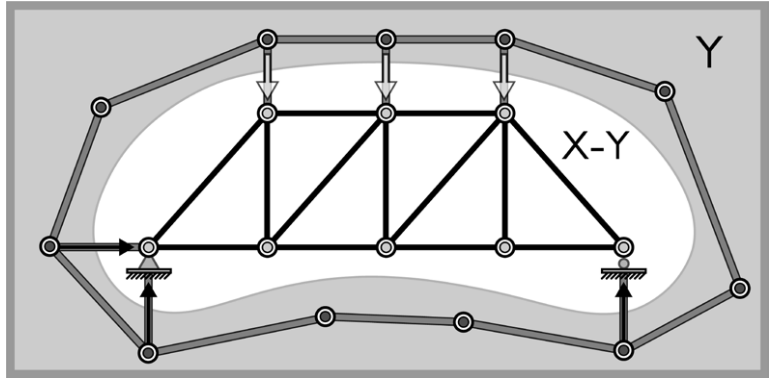

Figure 2. Boundary conditions in the form of loading and reaction forces are modeled as edges connecting the truss to an external loop $Y$.

### 5.3 Exact sequences

Exactness is a criterion of sequences of vector spaces and linear transformations that compactly characterizes a myriad of relationships between. Such a sequence of vector spaces

$$\ldots \to V_3 \to V_2 \to V_1 \to V_0 \to V_{-1} \ldots \tag{6}$$

is said to be *exact* if the image of every incoming transformation is equal to the kernel of the outgoing transformation at every vector space. That is, the homology of the sequence (6) completely vanishes. In the special case that (6) is non-zero only at the central three terms, there is an isomorphism of vector spaces $V_1 \cong V_0 \oplus V_2$ where $V_0 \cong V_1/V_2$; this can be worked out using properties of exactness [26].

It is a fact that any triple of cosheaves of the form $\mathcal{F} \hookrightarrow \mathcal{G} \to \mathcal{G}/\mathcal{F}$, where the first cosheaf map is an inclusion and the second is a quotient projection, induces a *long exact sequence* in homology

$$\ldots \to H_{k+1}\mathcal{G}/\mathcal{F} \to H_k\mathcal{F} \to H_k\mathcal{G} \to H_k\mathcal{G}/\mathcal{F} \to H_{k-1}\mathcal{F} \to \cdots, \tag{7}$$

derived from the corresponding exact sequence of cosheaf chain complexes [7]. True to its name, sequence (7) is an exact sequence of vector spaces and is a core feature of *homological algebra*.

Applying this to the externally loaded force cosheaf $\mathcal{F}_Y \hookrightarrow \mathcal{F}_X \to \mathcal{F}_{X-Y}$, where $\mathcal{F}_{X-Y} = \mathcal{F}_X/\mathcal{F}_Y$, one has

$$0 \to H_1\mathcal{F}_X \to H_1\mathcal{F}_{X-Y} \to H_0\mathcal{F}_Y \to H_0\mathcal{F}_X \to H_0\mathcal{F}_{X-Y} \to 0 \tag{8}$$

where $H_1\mathcal{F}_Y = 0$ (because $Y$ is assumed to be a loop with no self-stress) and all higher-homology terms vanish. This sequence collates the following elements:

1. The 1$^{\text{st}}$ homology $H_1\mathcal{F}_X$ is the vector space of *internal self*-stresses of $X$.

2. The relative 1st homology $H_1\mathcal{F}_{X-Y}$ is the total space of equilibrium stresses including those from boundary conditions.
3. The 0th homology $H_0\mathcal{F}_Y$ is the space of degrees of freedom of the frontier $Y$.
4. The 0th homology $H_0\mathcal{F}_X$ captures degrees of freedom of $X$ (potentially including those of $Y$).
5. The 0th homology $H_0\mathcal{F}_{X-Y}$ is the space of *constrained degrees of freedom* of $X - Y$ with respect to the various loadings and reactions restricting said freedom.

Exactness of (8) means that there are constraints tying together sequential terms. For example, the statement that *"every self-stress is an equilibrium stress with zero loading"* is equivalent to exactness of (8) at $H_1\mathcal{F}_X$ (the transformation $H_1\mathcal{F}_X \to H_1\mathcal{F}_{X-Y}$ is injective). Exactness at $H_1\mathcal{F}_{X-Y}$ has the interpretation: *"every equilibrium stress of a loaded structure is a sum of internal self-stresses and those resulting from outside forces"*. These outside forces are generated by instabilities in $Y$ that, in the larger structure $X$, "push" and "pull" on $X - Y$. Similar statements may be made about the remaining terms.

## 6. Planar graphic statics

Maxwell's seminal paper [14] made the astute observation that force vectors within members of a static truss combine to form a *reciprocal diagram* topologically dual to the primal (see Figure 3 for examples). Various aspects of 2D duality have been studied, such as variations on the planarity condition [21] or kinematic properties [19], [20]. Here we develop the root of the duality in its substantial form and extend to new territory.

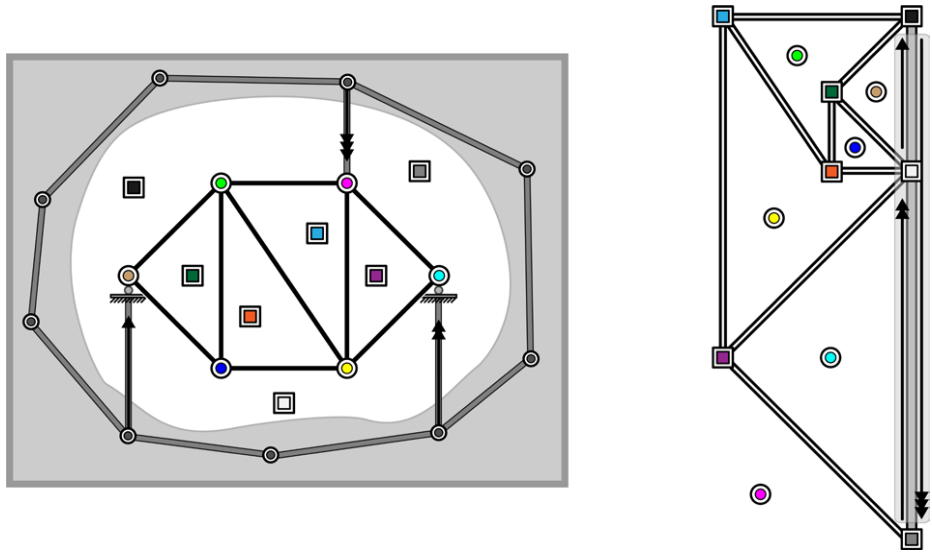

Figure 3. An example of graphic statics with boundary conditions is pictured. A form diagram is subject to external loads (left), with a single degree of equilibrium. A representative force diagram (a combined free body diagram) is drawn (right), with edge lengths equal to the load over the primal edge. Equilibrium of forces on the left is equivalent to dual face cells on the right being closed. The grey exterior (*cf.* Figure 2) on the left is dualized to the interior of the vertical degenerate triangle on the far right.

We focus on finite *planar graphs* with straight, non-overlapping edges. For a particular realization of a planar graph in the plane, the polyhedral regions interior to minimal cycles are *face cells*. Introducing these face cells extends a planar graph to a two-dimensional *cell complex*. Including the exterior which extends to infinity as one additional face yields a cell complex with the topology of the sphere $\mathbb{S}^2$. We remind the reader that although this cell complex $X = (V, E, F)$ has spherical topology, it remains realized in the flat plane by an assignment $p: V \to \mathbb{R}^2$. This is called a *form diagram*.

It is natural to extend our definitions of cosheaves over the plane as a 2D cell complex. The force cosheaf $\mathcal{F}$ over this spherical cell complex $X = (V, E, F)$ is an extension of that over its graph $(V, E)$: the stalks over all face cells are set to zero (0-dimensional vector spaces) with zero-matrices from face stalks to incident edge stalks. There is more to the story, based on how this force cosheaf is situated in the plane.

**6.1 The position cosheaf**

While the force cosheaf $\mathcal{F}$ describes the compatibility of forces over the truss, it is a natural question what the complementary data to forces are in the plane. Namely, since edge stalks $\mathcal{F}_e$ contain force information *parallel* to edges, what data structure describes information *perpendicular* to edges? This question naturally leads to its *complimentary cosheaf*, the quotient of the constant cosheaf $\mathbb{R}^2_X$ by $\mathcal{F}$. The reason for considering the constant cosheaf here is that it encompasses the *ambient space* of the plane. Stalks of $\mathbb{R}^2_X$ are akin to the tangent space: copies of $\mathbb{R}^2$ allocated to each location in the plane, stitched together by identity matrices.

There is an inclusion map $\phi: \mathcal{F} \hookrightarrow \mathbb{R}^2_X$ of cosheaves which we now describe. The force cosheaf has no edge stalks so $\phi_f$ is the zero matrix over faces. By $\phi_e$ over edges, the vector $\mathcal{F}_e \cong \mathbb{R}$ is sent to the vector $p(u) - p(v)$, the position of the edge endpoints. Over vertices, stalks of $\mathcal{F}$ and $\mathbb{R}^2_X$ are both copies of $\mathbb{R}^2$, so it is valid to let $\phi_v$ be the identity linear transformation. It is straightforward to check that $\phi$ is a valid cosheaf map by this construction.

The quotient cosheaf $\mathcal{G} \coloneqq \mathbb{R}^2_X / \mathcal{F}$, called the *position cosheaf*, is the complementary data structure we desire. Over faces, this cosheaf has face stalks $\mathcal{G}_f \cong \mathbb{R}^2$ encoding *positions* (to be detailed momentarily). There is trivial data over vertices, with $\mathcal{G}_v = 0$. The edge stalks $\mathcal{G}_e = \mathbb{R}^2 / \mathcal{F}_e$ represent the linear subspace perpendicular to the thrust of $e$ in the plane. The non-trivial cosheaf matrices $\mathcal{G}_{f \rhd e}$ are projection matrices $P_{e^\perp}$, projecting from the face stalk $\mathcal{G}_f$ onto the perpendicular space $\mathcal{G}_e$.

Recall the *Poincaré dual* cell complex to $X = (V, E, F)$ is the cell complex $\tilde{X} = (\tilde{F}, \tilde{E}, \tilde{V})$ where faces are converted to dual vertices, edges to dual edges, and vertices to dual faces [26]. In this way, the cell structure is dualized to an alternative cell structure on the sphere $\mathbb{S}^2$. We note that $\tilde{X}$ is merely the topological dual to $X$ as of now; the geometric bearing of cells of $\tilde{X}$ is informed by homology classes of the position cosheaf.

Face stalks of $\mathcal{G}$ encode the positions of the *dual vertices* of $\tilde{X}$ in $\tilde{F}$. Namely, $\mathcal{G}_f \cong \mathbb{R}^2$ may be considered as the space of possible positions for the dual vertex $\tilde{f}$ in the plane. With this identification, the homology space $H_2 \mathcal{G}$ consists of these realizations of dual vertices subject to constraints along edges. For every edge bounded by faces $e \lhd f_1, f_2$, there is the corresponding dual incidence $\tilde{e} \rhd \tilde{f_1}, \tilde{f_2}$ and dual restraint. For $\xi$ to be a homology element of $H_2 \mathcal{G}$ it must be the case that

$$\mathcal{G}_{f_1 \rhd e}(\xi_{f_1}) = \mathcal{G}_{f_2 \rhd e}(\xi_{f_2}) \tag{9}$$

or equivalently that

$$P_{e^\perp}(\xi_{f_1} - \xi_{f_2}) = 0. \tag{10}$$

In other words, the difference between $\xi_{f_1}$ and $\xi_{f_2}$ along the perpendicular direction to $e$ is zero – and so *any difference between these two positions lies along the linear subspace to edge $e$*. This means that the difference in positions of dual vertices $\tilde{f_1}, \tilde{f_2}$ is parallel to edge $e$ as determined by $\xi$. Moreover, every parallel dual realization $\xi: \tilde{F} \to \mathbb{R}^2$ is a homology class of $H_2 \mathcal{G}$ as this condition is descriptive. A parallel dual graph to that of Figure 1 is pictured in Figure 4.

In conclusion, the vector space $H_2 \mathcal{G}$ of the position cosheaf encodes *all possible parallel realizations* of the dual cell complex $\tilde{X}$ in the plane. These realizations of $\tilde{X}$ are called *force diagrams*.

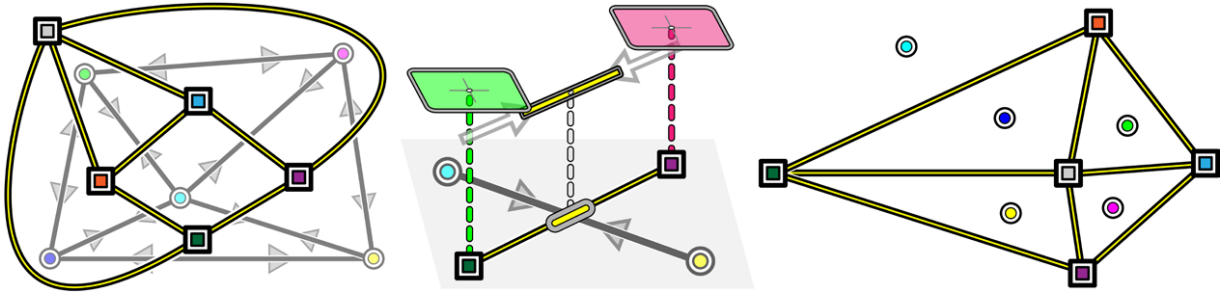

Figure 4. The topological dual to the form diagram of Figure 1 is pictured (left), with its geometry realized (right). A local view of the position cosheaf $\mathcal{G}$ is highlighted (center), with cosheaf matrices mapping face stalks to edge stalks. A generator of $H_2\mathcal{G}$ provides the geometric coordinates to dual vertices on the right, characterizing a force diagram. This parallel dual realization may be derived from the self-stress of the form diagram via the graphic statics of Theorem 1.

**6.2 Dual trusses**

We arrive at the heart of 2D graphic statics. While the classical history of form and force diagram duality reaches back over a century and a half [14], cosheaves permit an efficient modern reformulation. The cosheaf triple $\mathcal{F} \hookrightarrow \mathbb{R}_X^2 \to \mathcal{G}$ leads to a long exact sequence as per sequence (7):

$$\ldots \to H_2\mathcal{F} \to H_2\mathbb{R}_X^2 \to H_2\mathcal{G} \to H_1\mathcal{F} \to H_1\mathbb{R}_X^2 \to H_1\mathcal{G} \to H_0\mathcal{F} \to H_0\mathbb{R}_X^2 \to H_0\mathcal{G} \to 0 \qquad (11)$$

The force cosheaf has zero-stalks over faces, so $H_2\mathcal{F} \cong 0$. The constant cosheaf has homology $H_2\mathbb{R}_X^2 \cong \mathbb{R}^2$, $H_1\mathbb{R}_X^2 \cong 0$, and $H_0\mathbb{R}_X^2 \cong \mathbb{R}^2$, matching Betti numbers of the underlying sphere (multiplied by two, see example 1, Section 3.2). These computations simplify and split the exact sequence (11) into two:

$$0 \to \mathbb{R}^2 \to H_2\mathcal{G} \to H_1\mathcal{F} \to 0 \qquad (12)$$

$$0 \to H_1\mathcal{G} \to H_0\mathcal{F} \to \mathbb{R}^2 \to 0 \qquad (13)$$

The first exact sequence (18) *is* the establishment of 2D graphic statics. The space $\mathbb{R}^2$ corresponds to *global translations* of dual realizations of $\tilde{X}$ in $H_2\mathcal{G}$: these are global shifts of all dual vertices by the same translation vector. The linear transformation $H_2\mathcal{G} \to H_1\mathcal{F}$ takes the pointwise difference between dual coordinates and converts them to force values within a self-stress.

***Theorem 1:*** [14], [17] *Over a truss X realized in $\mathbb{R}^2$, there is a bijective correspondence between self-stresses of X and parallel dual realizations of $\tilde{X}$ in $\mathbb{R}^2$ up to global translation.*

This is the isomorphism of vector spaces

$$H_1\mathcal{F} \cong H_2\mathcal{G}/\mathbb{R}^2 \qquad (14)$$

following exact sequence (12). The ease of writing this equation conceals the underlying linear relationships at play. The force within each edge of a self-stress is equivalent to the distance between dual nodes in the corresponding force diagram.

A benefit of the abstraction presented here is the readily generalization of isolated graphic statics to graphic statics with boundary conditions – the isomorphism of vector spaces

$$H_1\mathcal{F}_{X-Y} \cong H_2\mathcal{G}_{X-Y}/\mathbb{R}^2. \qquad (15)$$

As $X - Y$ (Figure 2) has the topology of an open disk, its dual counterpart $\tilde{X} - \tilde{Y}$ has the topology of a closed disk. Equation (15) then is the statement that equilibrium stresses over a structure with boundary forces are equivalent to parallel realizations of the dual disk $\tilde{X} - \tilde{Y}$, up to global translation. An example is pictured in Figure 3.

Isomorphism (15) relies on $X - Y$ having $\beta_1 = 0$ so the constant cosheaf homology $H_1 \mathbb{R}_X^2$ vanishes. If this requirement is exempted it is still possible to formulate graphic statics, such as the case of *internally applied loadings* where loads are applied within the structure (and holes inserted into $X - Y$). In this case, exact sequence (11) does not split into the two sequences (12) and (13) so more care is required.

There are numerous applications and extensions of the abstract formulation here. With care, the boundary loop $Y$ can generalized to empower *relative graphic statics*. One can consider *partial* form diagrams $X - Y$ of a larger truss and rigorously compute consequent *partial* force diagrams $\tilde{X} - \tilde{Y}$. The relations to the whole are given by sequences analogous to (8).

**6.3 Impossible rotations of the dual**

The second exact sequence (13) also has meaning, although the homology $H_1 \mathcal{G}$ is difficult to interpret. This space consists of *impossible dual edge rotations* in the sense that edges of any parallel dual realization of $\tilde{X}$ cannot possibly rotate in the specified manner (even if edges are permitted to freely extend or contract axially). These edge rotations *obstruct* the creation of dual realizations. Just as mechanisms of a truss are not resistible by truss members, these rotations are *not realizable* by nodal positions. In fact, rotations of truss members under mechanism or global rotation *are* impossible edge rotations on dual force diagrams. Only the space $\mathbb{R}^2$ of global translations of $X$ does not rotate edges; the isomorphism

$$\mathbb{R}^2 \cong H_0 \mathcal{F} / H_1 \mathcal{G} \tag{16}$$

of vector spaces following exact sequence (13) certifies this. This relationship is novel.

**Theorem 2:** *The system of edge rotations from any mechanism or global rotation of a truss $X$ in $\mathbb{R}^2$ are impossible to satisfy on any parallel force diagram $\tilde{X}$ in $\mathbb{R}^2$, even while permitting dual edges in $\tilde{X}$ to freely extend axially.*

Note that dual edges have their rotations flipped when realized in tension versus compression (see Figure 5). The symmetry of sequences (12) and (13) suggests this theorem is as fundamental as Theorem 1. Indeed, if the values of an impossible edge rotation are converted to stress values along edges (with clockwise rotation converted to compression, and counterclockwise to tension for instance), they form a self-stress of the dual $\tilde{X}$. This is *graph reciprocity*, the mirroring of the force and form diagrams. This observation is in line with previous work involving mechanisms and duals diagrams [19], [20].

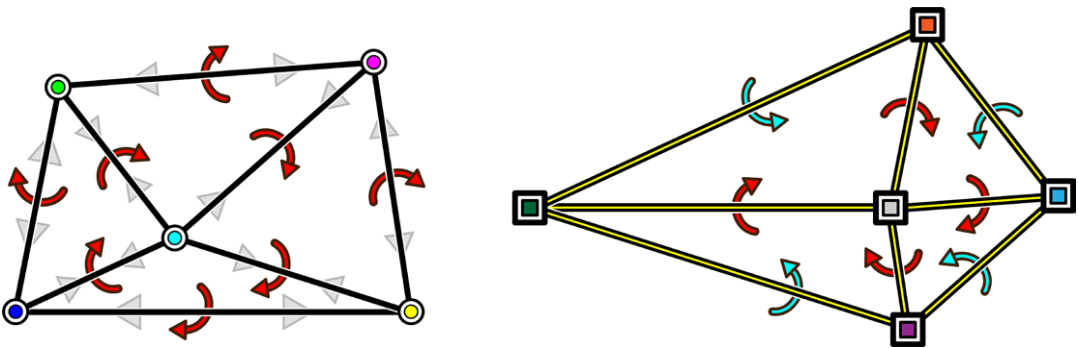

Figure 5. The form diagram of Figure 1 is globally rotated clockwise (left), inducing edge rotations. In a force diagram (right), dual edges have their rotations flipped when the primal edge is in tension. Even allowing dual edges to extend or contract, it is impossible to simultaneously rotate the internal four dual edges clockwise and rotate the external four dual edges counterclockwise. This illustrates Theorem 2; the single dimension of global rotation of the form diagram corresponds to the dimension of impossible dual edge rotations of force diagrams.

If the nodes of a truss are fixed connections instead of pinned, free-rotational connections, then impossible edge rotations are equivalent to states of *self-shear* of the structure. In this setting, perpendicular shear forces are applied along each member instead of parallel axial forces. These self-shear forces are in equilibrium, with zero net force at nodes.

One central application of graphic statics is in *form finding*, where a force diagram is manipulated to find efficient geometry for the primal truss. Sequence (13) provides a useful guarantee that changes to the force diagram do not activate instabilities of the form diagram. If $\xi$ is any change of coordinates of dual vertices of $\tilde{X}$ (an element of $C_2\mathcal{G}$), then $\partial\xi$ is in the zero class of $H_1\mathcal{G}$, thus zero in $H_0\mathcal{F}$.

**Theorem 3:** *For $X$ a form diagram and $\tilde{X}$ a force diagram embedded in $\mathbb{R}^2$, any infinitesimal repositioning of nodes of $\tilde{X}$ induces no translational, rotational, or mechanical motion of $X$.*

## 7. Beyond the plane

There are numerous potential extensions of the cosheaf-based graphic statics presented here.

One clear path is the extension of the well-known 3D polyhedral theory by Maxwell [14], [15] and Cremona [28], among others. Self-stresses of a cell complex $X$ in the plane are in bijective correspondence to a family of three-dimensional *lifted polyhedra* that project vertically onto $X$. As 2D reciprocal graphic statics is very much linked to plane liftings of the truss [17] with a correspondence that is homological in nature [18], one can build a cosheaf adaptation.

Recalling Example 2 on polynomial splines, there is potential to integrate splines with *linear shape functions* utilized in structural engineering. For a beam in bending, the deflection, moment, and shear quantities all may be approximated by polynomials. At ends, beams are connected and restrained by $C^r$ differentiability conditions - equivalent to compatibility conditions of spline cosheaves. More generally, the *finite element method* deconstructs complex shapes into simplified constituent components. It may be fruitful to investigate the homological properties of equilibrium equations in these models.

Finally, as presented in Sections 2 and 3, the force cosheaf $\mathcal{F}$ makes sense in dimension 3 and above. Indeed, the derivation of the Maxwell Rule (2) here given is no more difficult in $n$-D than in 2-D. To what extent, then, can the framework of cosheaves be integrated with recent breakthroughs in 3D polyhedral graphic statics [22]–[24]?

With 3-dimensional cells, the complexity between relationships of equilibrium spaces increases dramatically. While a single cosheaf map initiated 2D graphic statics in Section 6, higher dimensional graphic statics requires *filtrations* (or increasing sequences) of cosheaves, the analysis of which prompts the use of *spectral sequences*, a far-reaching generalization of the cosheaf long exact sequence (7).

## Acknowledgements

The authors appreciate insightful discussions with J. Hansen, M. Akbarzadeh, M. Hablicsek, E. Saliklis, and D. Lee. This work is supported by the Air Force Office of Scientific Research under award number FA9550-21-1-0334.

## References

[1] R. Ghrist, *Elementary Applied Topology ed 1.0*. Createspace, 2014.
[2] R. Ghrist, "Barcodes: the persistent topology of data," *Bulletin of the American Mathematical Society*, vol. 45, no. 1, pp. 61–75, 2008.
[3] A. E. Sizemore, J. E. Phillips-Cremins, R. Ghrist, and D. S. Bassett, "The importance of the whole: topological data analysis for the network neuroscientist," *Network Neuroscience*, vol. 3, no. 3, pp. 656–673, 2019.
[4] V. De Silva, R. Ghrist, and others, "Homological sensor networks," *Notices of the American mathematical society*, vol. 54, no. 1, 2007.

[5] M. Robinson, *Topological Signal Processing*, vol. 81. Springer, 2014.
[6] G. E. Bredon, *Sheaf Theory*. Springer New York, 1997.
[7] J. Curry, "Sheaves, Cosheaves and Applications," University of Pennsylvania, 2013.
[8] R. Ghrist and Y. Hiraoka, "Network codings and sheaf cohomology," *IEICE Proceedings Series*, vol. 45, no. A4L-C3, 2011.
[9] R. Ghrist and S. Krishnan, "Positive Alexander duality for pursuit and evasion," 2017.
[10] J. Hansen and R. Ghrist, "Distributed optimization with sheaf homological constraints," in *2019 57th annual allerton conference on communication, control, and computing (allerton)*, 2019, pp. 565–571.
[11] J. Hansen and R. Ghrist, "Opinion dynamics on discourse sheaves," *SIAM Journal on Applied Mathematics*, vol. 81, no. 5, pp. 2033–2060, 2021.
[12] H. Crapo, "Applications of geometric homology," in *French Workshop on Geometry and Robotics*, Springer, 1988, pp. 213–224. doi: 10.1007/3-540-51683-2_32.
[13] J. Hansen and R. Ghrist, "Toward a spectral theory of cellular sheaves," *Journal of Applied and Computational Topology*, vol. 3, no. 4, pp. 315–358, 2019.
[14] J. C. Maxwell, "XLV. On reciprocal figures and diagrams of forces," *The London, Edinburgh, and Dublin Philosophical Magazine and Journal of Science*, vol. 27, no. 182, pp. 250–261, 1864, doi: 10.1080/14786446408643663.
[15] J. C. Maxwell, "I.—On Reciprocal Figures, Frames, and Diagrams of Forces," *Transactions of the Royal Society of Edinburgh*, vol. 26, no. 1, pp. 1–40, 1870, doi: 10.1017/S0080456800026351.
[16] W. Whiteley, "Motions and Stresses of Projected Polyhedra," *Structural Topology*, vol. 7, no. June, pp. 13–38, 1982.
[17] H. Crapo and W. Whiteley, "Plane self stresses and projected polyhedra I : The Basic Pattern," *Structural Topology*, vol. 20, 1993.
[18] J. E. Hopcroft and P. J. Kahn, "A paradigm for robust geometric algorithms," *Algorithmica*, vol. 7, no. 1–6, pp. 339–380, 1992, doi: 10.1007/BF01758769.
[19] T. Mitchell, W. Baker, A. McRobie, and A. Mazurek, "Mechanisms and states of self-stress of planar trusses using graphic statics, part I: The fundamental theorem of linear algebra and the Airy stress function," *International Journal of Space Structures*, vol. 31, no. 2–4, pp. 85–101, 2016, doi: 10.1177/0266351116660790.
[20] A. McRobie, W. Baker, T. Mitchell, and M. Konstantatou, "Mechanisms and states of self-stress of planar trusses using graphic statics, part II: Applications and extensions," *International Journal of Space Structures*, vol. 31, no. 2–4, pp. 102–111, 2016, doi: 10.1177/0266351116660791.
[21] T. V. Mele and P. Block, "Algebraic graph statics," *CAD Computer Aided Design*, vol. 53, pp. 104–116, 2014, doi: 10.1016/j.cad.2014.04.004.
[22] M. Akbarzadeh, "3D Graphical Statics Using Reciprocal Polyhedral Diagrams," 2016. doi: 10.3929/ethz-a-010782581.
[23] M. Hablicsek, M. Akbarzadeh, and Y. Guo, "Algebraic 3D graphic statics: Reciprocal constructions," *CAD Computer Aided Design*, vol. 108, pp. 30–41, 2019, doi: 10.1016/j.cad.2018.08.003.
[24] M. Konstantatou, P. D'Acunto, and A. McRobie, "Polarities in structural analysis and design: n-dimensional graphic statics and structural transformations," *International Journal of Solids and Structures*, vol. 152, pp. 272–293, 2018.
[25] L. J. Billera, "Homology of smooth splines: generic triangulations and a conjecture of Strang," *Transactions of the american Mathematical Society*, vol. 310, no. 1, pp. 325–340, 1988.
[26] A. Hatcher, *Algebraic Topology*. Cambridge University Press, 2002.
[27] C. R. Calladine, "Buckminster Fuller's 'Tensegrity' structures and Clerk Maxwell's rules for the construction of stiff frames," *International Journal of Solids and Structures*, vol. 14, no. 2, pp. 161–172, 1978, doi: 10.1016/0020-7683(78)90052-5.
[28] L. Cremona, *Le figure reciproche nella statica grafica*. Tipografia di G. Bernardoni, 1872.